\documentclass{amsproc}

\usepackage{amsthm,amsfonts}
\usepackage{amsmath,amssymb,amscd}
\DeclareMathOperator{\Coh}{Coh}
\DeclareMathOperator{\coker}{coker}
\DeclareMathOperator{\End}{End}
\DeclareMathOperator{\Ext}{Ext}
\DeclareMathOperator{\GL}{GL}
\DeclareMathOperator{\hd}{hd}
\DeclareMathOperator{\Hom}{Hom}
\DeclareMathOperator{\HOM}{\mathcal{H}om}
\DeclareMathOperator{\ob}{ob}

\newtheorem{thm}{Theorem}[section]
\newtheorem{lemma}[thm]{Lemma}
\newtheorem{prop}[thm]{Proposition}
\newtheorem{cor}[thm]{Corollary}
\theoremstyle{definition}

\theoremstyle{remark}
\newtheorem{rem}[thm]{Remark}

\numberwithin{equation}{section}
\begin{document}

\title{Homological Mirror Symmetry in Dimension One}

\author{Bernd Kreu{\ss}ler}
\address{FB Mathematik, Universit\"at Kaiserslautern, 67653 Kaiserslautern,
  Germany} 
\email{kreusler@mathematik.uni-kl.de}

\subjclass{14J32, 14H52, 53D12} 

\begin{abstract}
In this paper we complete the proof began by A. Polishchuk and E. Zaslow
\cite{PZ} of a weak version of Kontsevich's homological mirror symmetry
conjecture for elliptic curves.  
\end{abstract}

\maketitle

\section{Introduction}\label{sec:intro}

This note grew out of an attempt to understand the details in the paper of
A. Polish\-chuk and E. Zaslow \cite{PZ} on the homological mirror symmetry
conjecture of M. Kontsevich. This conjecture was formulated in his ICM talk
\cite{KonICM} in Z\"urich in 1994. Roughly speaking, this conjecture describes
a close relation between two seemingly totally unrelated mathematical
structures defined in dependence of two Calabi-Yau manifolds, a so-called
mirror pair. Although the conjecture is mostly interesting in dimension three,
it is mathematically interesting also in other dimensions. 

Of course, one  encounters the easiest situation in dimension one. A one
dimensional Calabi-Yau manifold is just a two dimensional torus (elliptic
curve). In \cite{PZ} precisely this case is investigated. There exist attempts
to generalize their results to higher dimensional tori, undertaken by C. van
Enckevort \cite{Enc} and K. Fukaya \cite{Fuk}. These papers show that the
conjecture is very hard, even in the presumably simple case of a torus. To date
there is no complete proof of the conjecture, even in dimension one. In
\cite{PZ} a weaker modification of Kontsevich's conjecture is studied, but it
is not completely clear how to obtain a similar result in the higher
dimensional case without restricting additionally the objects in the categories
compared. 

The contribution of this paper to the conjecture is a complete proof and a
correct formulation of the weakened version of Kontsevich's conjecture
described in \cite{PZ}. We will not repeat all steps of the proof. Here we
restrict to those details whose proofs seem to require more explanation or
are lacking. The author hopes to contribute in this way to a better
understanding of the beautiful results of \cite{PZ}. The availability of
\cite{Enc} and \cite{Fuk} was particularly helpful in preparing this paper.

The homological mirror symmetry conjecture, as stated by M. Kontsevich in his
ICM talk in Z\"urich in 1994 \cite{KonICM} can be formulated as follows:

For every Calabi-Yau (3-)manifold $X$ there exists a mirror partner $X^\circ$
with a symplectic form $\omega$ and an equivalence of ($A_\infty$)-categories:
$$\mathcal{D}^b(X) \cong \mathcal{F}(X^\circ).$$

We denote by $\mathcal{D}^b(X)$ the derived category of the abelian category of
coherent sheaves on $X$ (\cite{GM}, \cite{HG}).  By $\mathcal{F}(X^\circ)$ we
denote Fukaya's category, which actually is not a category in the usual
sense. It is a so-called $A_\infty$-category. To give the conjecture a precise
meaning, one has to either equip $\mathcal{D}^b(X)$ with an
$A_\infty$-structure or to modify Fukaya's category in such a way that it
becomes a triangulated (ordinary) category. In this paper we follow \cite{PZ}
and replace $\mathcal{F}(X^\circ)$ by its cohomology and compare it with
$\mathcal{D}^b(X)$ as an additive category (forgetting the triangulated
structure). To prove deeper versions of mirror symmetry in the one dimensional
case involving the $A_\infty$-structure, it will probably be helpful to have a
detailed proof in this simple case. A proof of the equivalence of two
$A_\infty$-structures on the category of vector bundles on an elliptic curve
was given in \cite{Pinfty}. One of these structures is induced via the
equivalence with the cohomology of Fukaya's category and the second
$A_\infty$-structure was introduced in \cite{Phigh}, inspired by
\cite{Merk}. But in \cite{Pinfty} only so-called transversal structures are
considered, non-transversal Lagrangians are not taken into account.

In Section \ref{sec:derived} we collect basic facts about the structure of the
bounded derived category of coherent sheaves on a smooth projective curve
necessary for the proof of the main theorem. Because the construction of our
equivalence of categories relies heavily on the use of a version of Serre
duality for coherent sheaves on an elliptic curve, a proof of it is presented
here. 

In Section \ref{sec:fukaya} we recall the definition of Fukaya's category as
given in \cite{PZ}. To overcome the fact that \cite{PZ} do not define spaces of
morphisms between arbitrary pairs of objects, we define the cohomology
$\mathcal{F}^0(X)$ of Fukaya's category directly, i.e. without giving a
definition of Fukaya's $A_\infty$-category $\mathcal{F}(X)$. To be able to
prove the main theorem we need to enlarge this category to get an additive
category $\mathcal{FK}^0(X)$. We study some functors between such
categories and collect important properties of them in analogy to the
statements in Section \ref{sec:derived}.

Section \ref{sec:equiv} is devoted to the proof of the main theorem stating
an equivalence of categories 
$$\mathcal{D}^b(X) \cong \mathcal{FK}^0(X^\circ).$$ 
We focus on those parts of the proof which are not contained in the paper of
A. Polishchuk and E. Zaslow to keep the paper short. 

\textsl{Acknowledgement:} The author would like to thank all the people at the
Department of Mathematics at Oklahoma State University for their hospitality
during his one-year visit in Stillwater, where the main work on this paper
was done. Particular thanks go to Sheldon Katz for stimulating discussions and
generous hospitality. 

\section{The derived category}
\label{sec:derived}

One reason, that makes possible the proof of a weak version of the homological
mirror symmetry conjecture in dimension one, is the fact that the bounded
derived category $\mathcal{D}^b(X):= \mathcal{D}^b(\Coh_X)$ of a smooth
projective curve $X$ has a well understood structure. To become more specific
let us recall \cite{KS} that the homological dimension $\hd{\mathcal{A}}$ of an
abelian category $\mathcal{A}$ is by definition the smallest non-negative
integer $k$ such that for all $j>k$ and all objects $A, B$ in
$\mathcal{A}:\quad \Ext^j(A,B)=0.$

We know, for example, that the category $\mathcal{A}=\Coh_X$ of coherent
sheaves on a smooth projective variety $X$ has homological dimension
$\hd(\Coh_X) = \dim{X}$, because any coherent sheaf on $X$ has a resolution 
of length at most equal to $\dim{X}$ whose members are vector bundles of finite
rank.

On the other hand, the following lemma, which is an easy exercise (see
\cite{KS}), provides the information we are interested in. Here we denote as
usual by $A[-k]$ the object in $\mathcal{D}^b(\mathcal{A})$ given by the
complex having the object $A$ of the abelian category $\mathcal{A}$ at position
$k$ and the zero object elsewhere.

\begin{lemma}
  If $\mathcal{A}$ is an abelian category with $\hd(\mathcal{A})\le
  1$ and $A^{\bullet}$ an object of $\mathcal{D}^b(\mathcal{A})$, then
  $A\cong\bigoplus_k H^k(A^{\bullet})[-k]$ in $\mathcal{D}^b(\mathcal{A})$.
\end{lemma}

This gives the following structure result for the bounded derived category of a
curve:

\begin{cor}
  If $X$ is a smooth projective curve and $A^{\bullet}$ an object of
  $\mathcal{D}^b(X)$, then $A\cong\bigoplus_k   H^k(A^{\bullet})[-k]$ in
  $\mathcal{D}^b(X)$. 
\end{cor}

Here we denote by $H^k(A^{\bullet})$ the $k$-th cohomology sheaf of the complex
$A^{\bullet}$ (considered as an object in the derived category). It is not the
sheaf cohomology of the members of this complex. 

Call a coherent sheaf indecomposable if it is not the direct sum of two
non-trivial coherent sheaves. Now we can formulate the following result.

\begin{prop}
  Let $X$ be a smooth projective curve. By taking finite direct sums of objects
  of the form $A[n]$, where $A$ is an indecomposable coherent sheaf on $X$, we
  obtain all objects of $\mathcal{D}^b(X)$ up to isomorphism. 
\end{prop}

By using standard results from category theory \cite{ML} we can reformulate
this in the following way:

\begin{prop}\label{subcat}
  Let $X$ be a smooth projective curve. The full subcategory of
  $\mathcal{D}^b(X)$ formed by all objects which are finite direct sums of
  objects of the form $A[n]$, where $A$ is an indecomposable coherent sheaf on
  $X$, is equivalent to $\mathcal{D}^b(X)$. 
\end{prop}

We will work with this subcategory instead of the bounded derived category
itself, but for simplicity we will not always mention this explicitly. Note
this subcategory is closed under shifts and finite direct sums.

On a smooth projective curve $X$ for any coherent sheaf $A$ the torsion free
sheaf $A/A_{\text{tor}}$ is locally free and because a coherent sheaf with zero
dimensional support has vanishing first cohomology, any coherent sheaf is the
direct sum of a locally free sheaf and a torsion sheaf: $A\cong
A/A_{\text{tor}}\oplus A_{\text{tor}}$. In particular, any indecomposable
coherent sheaf is either an (indecomposable) torsion sheaf supported at one
point or an (indecomposable) vector bundle. This observation will be used
frequently.

We need a more explicit description of the indecomposable coherent sheaves on
elliptic curves. We denote by $E_\tau$ the elliptic curve defined by the
lattice $\Gamma_\tau=\mathbb{Z}\oplus\tau\mathbb{Z}$ where $\tau$ is a complex
number lying in the upper half plane. By abuse of language we will speak about
``the point $a\tau+b$'' on $E_\tau$ if we mean the point $a\tau+b \mod
\Gamma_\tau$. 

All line bundles on $E_\tau$ are obtained in the following way: Let $\varphi:
\mathbb{C}\rightarrow \mathbb{C}^\ast$ be a holomorphic function satisfying
$\varphi(z+1) = \varphi(z)$ for all $z\in\mathbb{C}$. The line
bundle $\mathcal{L}(\varphi)$ is by definition
$(\mathbb{C}\oplus\mathbb{C})/\Gamma_\tau$ 
where the action of the lattice is defined by $\tau\cdot(z,v) = (z+\tau,
\varphi(z)v)$ (and $1\in\Gamma_\tau$ acts trivially). If it is important to
stress the dependence on $\tau$ we write $\mathcal{L}_\tau(\varphi)$ instead of
$\mathcal{L}(\varphi)$. A special role is played by the function $\varphi_0(z)
= \exp(-\pi i\tau-2\pi iz)$. The line bundle $\mathcal{L}_\tau(\varphi_0)$
corresponds to the divisor of degree one supported at the point
$\frac{1+\tau}{2}$. We obtain all line bundles on $E_\tau$ by considering
functions $\varphi(z)=t^\ast_x\varphi_0^{\phantom{1}}\cdot\varphi_0^{n-1}$,
where $x=a\tau+b$ is a point on $E_\tau$, $t_x$ is the translation by $x$
(i.e. $t^\ast_x\varphi(z) = \varphi(z+x)$) and $n$ is an integer. The global
sections of $\mathcal{L}(\varphi)$ are given by holomorphic functions
$f:\mathbb{C}\rightarrow\mathbb{C}$ fulfilling the identities $f(z+\tau) =
\varphi(z)f(z)$ and $f(z+1)=f(z)$ for all $z\in \mathbb{C}$. A basis of the
space of global sections
$H^0(\mathcal{L}(t^\ast_x\varphi_0^{n}))$ is formed by the
theta functions with characteristics
$t^\ast_x\theta\big[\frac{k}{n},0\big](n\tau,nz)$ where $0\le k < n$. In
particular, the dimension of this space of sections is $n$ if $n> 0$. According
to Mumford \cite{Mum} we use here 
$$\theta\big[a,b\big](\tau,z) := \sum_{n\in\mathbb{Z}} \exp(\pi i(n+a)^2 + 2\pi
i(n+a)(z+b)).$$

To describe vector bundles of higher rank we use the notion of a unipotent
vector bundle. A vector bundle $F$ is called \emph{unipotent} if it has a
filtration $0=F_0\subset F_1\subset \ldots\subset F_r=F$ with $F_{i+1}/F_i\cong
\mathcal{O}_{E_\tau}$. We will consider the following unipotent bundles: Let
$V$ 
be a (finite dimensional) complex vector space and $A\in\GL(V)$, then the
bundle $F(V,A)$ is by definition $(\mathbb{C}\times V)/\Gamma_\tau$ where the
action of $\Gamma$ is given by $\tau\cdot(z,v) = (z+\tau, Av)$ (and
$1\in\Gamma$ acting trivially). In the particular case where $A=\exp(N)$ with a
nilpotent endomorphism $N\in\End(V)$ we obtain a unipotent bundle $F(V,\exp(N))
= F_\tau(V,\exp(N))$. If the kernel of $N$ is of dimension one, then
$F(V,\exp(N))$ is indecomposable. A nilpotent endomorphism with one dimensional
kernel will be called \emph{cyclic}, because it defines the structure of a
cyclic $\mathbb{C}[T]$-module on $V$. 

Later we shall need the following lemma, which follows easily from the fact
that holomorphic sections of a flat vector bundle are covariantly constant.

\begin{lemma}\label{flatsect}
  If $A\in\GL(V)$ then $$H^0(E_\tau,F_\tau(V,A)) = \ker(\mathbf{1}_V-A).$$ In
  particular, $$\Hom(F(V_1,A_1),F(V_2,A_2)) = \{f\in\Hom(V_1,V_2)\mid f\circ
  A_1=A_2\circ f\}.$$
\end{lemma}

Thanks to the following result, which follows from work of M.F. Atiyah
\cite{At}, this is enough to have an explicit description of all indecomposable
vector bundles on an elliptic curve. To formulate the result we need the notion
of an isogeny. This is a surjective homomorphism of algebraic groups with
finite kernel. We need only the isogenies:
$\pi_r:E_{r\tau}\rightarrow E_\tau$ coming from the inclusion
$\Gamma_{r\tau}\subset \Gamma_\tau$ and the identity map on $\mathbb{C}$.

\begin{thm}\label{vb} 
  Any indecomposable vector bundle on an elliptic curve $E_\tau$  is
  of the form $\pi_{r\ast}( \mathcal{L}_{r\tau}(\varphi)\otimes
  F_{r\tau}(V,\exp(N)))$, with $N$ a cyclic nilpotent endomorphism.  
\end{thm}

An indecomposable coherent sheaf which is not a vector bundle is a torsion
sheaf supported at one point $x=a\tau+b$. The space  $V$ of global sections of
such a sheaf is a finite dimensional complex vector space with the structure of
an $\mathcal{O}_{E_\tau,x}$-module. Since we are on a smooth curve, we can
choose a generator of the maximal ideal of $\mathcal{O}_{E_\tau,x}$ (a
uniformizing parameter). The module structure is then determined by the action
of this generator on $V$. Because we constructed $E_\tau$ as a quotient of
$\mathbb{C}$ we can use the function $z-x$ as a generator, where $z$ denotes
the usual holomorphic coordinate on $\mathbb{C}$. Thus, a coherent
sheaf supported at one point is determined by the point $x$ supporting it, a
vector space $V$ and a nilpotent endomorphism $N\in\End(V)$. We denote it by
$S(x,V,N)$ or sometimes also by $S_\tau(x,V,N)$. Again, the sheaf is
indecomposable if and only if $N$ has one dimensional kernel.

Because the author was not able to find the following version of Serre duality
in the literature, we give a proof here.

\begin{lemma}\label{ellserre}
  Let $A_1, A_2$ be two coherent sheaves on the elliptic curve $E_\tau$, then
  we have a functorial isomorphism
  $$\Ext^1(A_2,A_1)\cong\Hom(A_1,A_2)^\ast.$$
\end{lemma}

\begin{proof}
Recall that the dualizing sheaf of a smooth elliptic curve is
the structure sheaf. Using Serre duality as presented in \cite[III.7.6]{H} and
using \cite[III.6.3, III.6.7]{H} we obtain the statement of the Lemma with the
additional assumption that $A_2$ is locally free. 

Let now $A_2$ be an arbitrary coherent sheaf. Because $E_\tau$ is a smooth
projective curve, there exists an exact sequence $0\rightarrow F'\rightarrow
F\rightarrow A_2\rightarrow 0$ with locally free sheaves $F$ and $F'$. By
applying the covariant functors $\Ext^1(A_1,\cdot)$ and
$\Hom(\cdot,A_1)^\ast$ we obtain the following diagram with exact rows:

$$
\begin{CD}
  \Ext^1(A_1,F') @>>> \Ext^1(A_1,F) @>>> \Ext^1(A_1,A_2) @>>> 0\\
  @V{\cong}VV @V{\cong}VV @V{\exists}VV\\
  \Hom(F',A_1)^\ast @>>> \Hom(F,A_1)^\ast @>>> \Hom(A_2,A_1)^\ast @>>> 0\\
  \end{CD}
$$

The commutativity of the left square ensures the existence of the right
vertical arrow, which is an isomorphism by the five lemma. By construction,
this isomorphism is functorial with respect to $A_1$. To get functoriality with
respect to $A_2$ we have to show commutativity of the diagram
\begin{equation}
  \label{func}
\begin{CD}
  \Ext^1(A_1,A_2) @>>> \Ext^1(A_1,\tilde{A}_2)\\
  @VV{\cong}V @VV{\cong}V \\
  \Hom(A_2,A_1)^\ast @>>> \Hom(\tilde{A}_2,A_1)^\ast\\
\end{CD}
\end{equation}
arising from a morphism of coherent shaves $A_2\rightarrow\tilde{A}_2$ and
resolutions $$0\rightarrow F'\rightarrow F\rightarrow A_2 \rightarrow 0$$ and 
$$0\rightarrow \tilde{F}'\rightarrow \tilde{F}\rightarrow \tilde{A}_2
\rightarrow 0$$ as above. A priori the vertical isomorphism in the diagram
(\ref{func}) might depend on the resolutions chosen. But once we have shown
commutativity for the identity $A_2\rightarrow A_2$ (but arbitrary resolutions)
we know that the isomorphism $\Ext^1(A_1,A_2) \cong \Hom(A_2,A_1)^\ast$ does
not depend on the resolution of $A_2$. In particular, if $A_2$ is locally free,
this is already known. 

To show commutativity of (\ref{func}) one usually lifts
$A_2\rightarrow\tilde{A}_2$ to a morphism of complexes:
$$
\begin{CD}
  0 @>>> F' @>>> F @>>> A_2 @>>> 0\\
  &&@VVV @VVV @VVV\\
  0 @>>> \tilde{F}' @>>> \tilde{F} @>>> \tilde{A}_2 @>>> 0\\
\end{CD}
$$
but this will not be possible in general. 

We consider first the situation where $A_2=\tilde{A}_2$ is a torsion sheaf,
$\tilde{F}=F\otimes L$ and $\tilde{F}'=F'\otimes L$ for a line bundle $L$ on
$E_\tau$. Observe that $A_2\otimes L\cong A_2$ since $A_2$ is a torsion
sheaf. If the degree of $L$ is large enough we have $\Ext^1(F,F'\otimes L) =
H^1(F^\vee\otimes F'\otimes L) = 0$ and obtain a lift for any morphism
$A_2\rightarrow A_2$. This shows that we obtain the same Serre duality
isomorphism, if we twist the resolution of a torsion sheaf by a line bundle of
sufficiently high degree.

If now $\tilde{A}_2$ is an arbitrary torsion sheaf we can apply the same
argument to show that we can lift any morphism $A_2\rightarrow\tilde{A}_2$ to a
morphism from $$0\rightarrow F'\rightarrow F\rightarrow A_2 \rightarrow 0$$ to
a twisted resolution $$0\rightarrow \tilde{F}'\otimes L \rightarrow
\tilde{F}\otimes L \rightarrow \tilde{A}_2 \rightarrow 0$$ where $L$ is of
sufficiently high degree. This shows commutativity of diagram (\ref{func}) in
case $\tilde{A}_2$ is torsion and so also the independence on the resolution of
the Serre duality isomorphism if $A_2$ is a torsion sheaf. 

If $\tilde{A}_2$ is locally free, any torsion direct summand of $A_2$ will be
mapped to zero in $\tilde{A}_2$. So we can assume $A_2$ to be locally free in
this case. Then the known functoriality of the Serre duality isomorphism
implies the commutativity of (\ref{func}). Using the additivity of the
functors involved, we obtain the required commutativity in general. 
\end{proof}

To understand morphisms and their composition in $\mathcal{D}^b(E_\tau)$ we
recall the functorial isomorphism 
$$\Hom_{\mathcal{D}^b(E_\tau)}(A_1[m], A_2[n]) =
\Hom_{\mathcal{D}^b(E_\tau)}(A_1, A_2[n-m]) = \Ext^{n-m}(A_1, A_2)$$ 
for objects $A_1, A_2$ in $\Coh_{E_\tau}$. Since $E_\tau$ is a curve this
vanishes if $n-m\notin\{0,1\}$. 

So we have $\Hom_{\mathcal{D}^b(E_\tau)}(A_1, A_2) = \Hom(A_1,A_2)$, which
means that $\Coh_{E_\tau}$ is a full subcategory of $\mathcal{D}^b(E_\tau)$ if
we send $A$ to $A[0]$. On the other hand, if we combine the above isomorphism
with Serre duality (Lemma \ref{ellserre}) we obtain a functorial isomorphism
$$\Hom_{\mathcal{D}^b(E_\tau)}(A_1, A_2[1]) = \Ext^1(A_1, A_2) = \Hom(A_2,
A_1)^\ast.$$ using functoriality in the first argument we obtain a commutative
diagram 
$$
\begin{CD}
  \Hom_{\mathcal{D}^b(E_\tau)}(A_1, A_2) \times
  \Hom_{\mathcal{D}^b(E_\tau)}(A_2, A_3[1]) @>{\circ}>>
  \Hom_{\mathcal{D}^b(E_\tau)}(A_1, A_3[1])\\
  @VV=V @VV=V\\
  \Hom(A_1,A_2) \times \Ext^1(A_2,A_3) @>>> \Ext^1(A_1,A_3)\\
  @VV{\cong}V @VV{\cong}V\\
  \Hom(A_1,A_2) \times \Hom(A_3, A_2)^\ast @>>> \Hom(A_3, A_1)^\ast\\
\end{CD}
$$
where the arrow in the bottom row sends $(f,\psi)$ to $\psi(f\circ\cdot) \in
\Hom(A_3,A_1)^\ast$. The map in the second row is usually called Yoneda
pairing. Similarly, using functoriality in the second argument, we translate
the composition $$\Hom_{\mathcal{D}^b(E_\tau)}(A_1, A_2[1]) \times
  \Hom_{\mathcal{D}^b(E_\tau)}(A_2[1], A_3[1]) \rightarrow
  \Hom_{\mathcal{D}^b(E_\tau)}(A_1, A_3[1])$$
to the map $$\Hom(A_2,A_1)^\ast \times \Hom(A_2,A_3) \rightarrow
\Hom(A_3,A_1)^\ast$$ sending $(\varphi, g)$ to $\varphi(\cdot\circ g) \in
\Hom(A_3, A_1)^\ast$.

\section{A version of Fukaya's category}
\label{sec:fukaya}

In his formulation of homological mirror symmetry Kontsevich used a modified
version of Fukaya's $A_\infty$-category. In Fukaya's original construction
\cite{FukMH} the objects were certain Lagrangian submanifolds of a fixed
symplectic manifold.  The morphisms in this $A_\infty$-category are defined
with the aid of Floer complexes for Lagrangian intersections. Inspired by some
developments in physics (D-branes), Kontsevich proposed to take as objects
special Lagrangian submanifolds equipped with unitary local systems. In the
paper \cite{Fuk}, Fukaya considers Hamiltonian isotopy classes of Lagrangian
submanifolds equipped with flat line bundles. At the moment it doesn't seem to
be clear which one is the final definition. We will follow here the paper
\cite{PZ} and use as objects pairs consisting of a special Lagrangian
submanifold and a local system on it whose monodromy has eigenvalues of modulus
one. It is an interesting question to study the relations between all these
different versions of Fukaya's category, even in the case of a torus.

There are some attempts to define an $A_\infty$-category for any symplectic
manifold (equipped with a $B$-field) \cite{FukMH}, \cite{Fuk}. But there are
some unclear technical details in the general case. In the case of a torus of
dimension one the situation is simple enough to be able to study the cohomology
of this $A_\infty$-category. In particular, the graded spaces of morphisms
between transversal Lagrangians consist of just one graded piece and so the
differential on this complex is trivial. Taking cohomology is then the same as
taking the zeroth graded piece. The appearance of this simple structure seems
to be a second reason for the accessibility of the proof. 

We are not treating
the $A_\infty$-structure here, because the space of homomorphisms between
objects with the same underlying Lagrangian submanifold has more that one
graded piece and is determined up to homotopy only. Results about Floer
homology suggest (see \cite{Fuk}) what the correct cohomology space of such a
complex should be and so we are able to define $\mathcal{F}^0(E^\tau)$ without
having $\mathcal{F}(E^\tau)$. See also \cite{Enc}.

We consider the two dimensional torus $\mathbb{T}^2=\mathbb{R}^2/\mathbb{Z}^2$
and equip it with a symplectic structure by using the volume form
$A(dx^2+dy^2)$, where $x$ and $y$ are real coordinates and $A>0$ is a positive
real number. The K\"ahler form of this metric is $Adx\wedge dy$. We choose, in
addition, a so-called $B$-field which is an element of $H^2(\mathbb{T}^2,
\mathbb{R})/H^2(\mathbb{T}^2, \mathbb{Z})$. It can be represented by a two form
$B dx\wedge dy$ with $B\in\mathbb{R}$. Choosing such a representative we obtain
a complex number $\tau=B+iA$ lying in the upper half plane. The $B$-field
determines $\tau$ modulo $\mathbb{Z}$. The torus $\mathbb{T}^2$ equipped with
$\tau$ will be denoted by $E^\tau$. The two form $\tau(dx\wedge dy)$ is also
called the \emph{complexified K\"ahler form}.

Let us recall the definition of $\mathcal{F}^0(E^\tau)$, the cohomology 
of Fukaya's $A_\infty$-category $\mathcal{F}(E^\tau)$ for the symplectic
torus. We will not consider $A_\infty$-structures here. We follow \cite{PZ} to
define the objects and morphisms in the transversal case. In contrast with the
definition of Polishchuk and Zaslow we define morphisms between any pair of
objects, not only in the transversal case.

The objects are triples $(\Lambda,\alpha,M)$ where
$\Lambda\subset E^\tau$ is a closed submanifold given by an affine line in
$\mathbb{R}^2$. We describe such a submanifold usually by a parametrization
$(x(t),y(t))$, where the two components depend linearly on the real parameter
$t$. To get a closed submanifold, this line must have rational slope. 

Such submanifolds equipped with an orientation are precisely the special
Lagrangian submanifolds of $E^\tau$. For a discussion of the concept of special
Lagrangian submanifolds see e.g. \cite{HL}.

To get a shift functor in Fukaya's category we have to equip the Lagrangian
submanifolds with an additional structure called a grading. A detailed study of
this concept can be found in \cite{Sei}. In general, this additional structure
is necessary to define the Maslov index, which is used to introduce a grading
on the spaces of morphisms in Fukaya's $A_\infty$-category. 

Because we consider special Lagrangian submanifolds of a two torus, things
simplify a lot. The additional structure is given by the choice of a logarithm
$\alpha$ of the slope of $\Lambda$. More precisely, $\alpha$ is a real number
with the property that there exists $t\in\mathbb{R}$ with
$x(0)+iy(0)+\exp(\pi i\alpha) = x(t)+iy(t)$. 

Finally, $M$ stands for a local system on $\Lambda$
whose monodromy has eigenvalues with unit modulus. This is slightly more
general than considering only unitary local systems. 

By a local system we mean here a locally constant sheaf of complex vector
spaces. A local system can be described equivalently by a representation of the
fundamental group of the underlying manifold or by a complex vector bundle
equipped with a flat connection. Because $\Lambda=S^1$ is a circle, its
fundamental group is a free abelian group with one generator. If we fix an
orientation on $\Lambda$, a generator of the fundamental group is determined
and a local system on $\Lambda$ is given by a vector space $V$ and an
automorphism $M\in\GL(V)$. If we change the orientation and consider $M^{-1}$,
we obtain an isomorphic local system on $\Lambda$. Two conjugate automorphisms
define isomorphic local systems. The operator $M$ is called the monodromy of
the local system and we are considering here only monodromy automorphisms whose
eigenvalues have modulus one. 

Considering a triple $(\Lambda,\alpha,M)$ and speaking about the local system
$M$ we mean the local system on $\Lambda$ determined
by $M\in\GL(V)$ and the orientation on $\Lambda$ given by the ``direction''
$\exp(\pi i\beta)$, where $\beta$ is the unique real number with
$-\frac{1}{2}<\beta\le\frac{1}{2}$ and $\beta-\alpha\in\mathbb{Z}$.

Sometimes it is convenient to replace the stalk $M_x$ of a local system $M$ by
the vector space $V$. To do this in a coherent way we introduce the following
convention: If $\Lambda=S^1$ is equipped with an orientation, the local system
given by $M\in\GL(V)$ can be described as the sheaf of (local) solutions $s$ of
the differential equation $\nabla s=0$ on the vector bundle
$\tilde{\Lambda}\times V/\sim$. Here we let $\tilde{\Lambda}\rightarrow\Lambda$
be the universal cover, which is the affine line in $\mathbb{R}^2$ defining
$\Lambda\subset E^\tau$ in our special situation. By $\lambda\in\mathbb{R}^2$
we denote the unique vector parallel to $\tilde{\Lambda}$ representing a lift
of the positive generator of $\pi_1(\Lambda)$. The equivalence relation $\sim$
is generated by $(l,v)\sim(l+\lambda, Mv)$. If $U\subset\Lambda$ is open and
$s:U\rightarrow V$ defines a local section of $\tilde{\Lambda}\times V/\sim$,
we define $\nabla s=ds$. This means, the local system is formed by locally
constant maps $U\rightarrow V$. If $x\in\Lambda$ and we choose
$\tilde{x}\in\tilde{\Lambda}$ over $x$ we get from this construction an
isomorphism $V\stackrel{\sim}{\rightarrow} M_x$ between $V$ and the stalk of
the local system $M$ at $x$. Our convention will be to choose once an for all
an interval of the form $\{\tilde{x}_0+t\lambda\mid 0\le t\le 1\}
\subset\tilde{\Lambda}$ for any $\Lambda$ and use the isomorphism $V
\stackrel{\sim}{\rightarrow} M_x$ defined by the point in this particular
interval which lies over $x$. The advantage of this choice is that in case
$M=\exp(M')$ the monodromy of the local system $M$ along a path
$\gamma:[0,1]\rightarrow \Lambda$ is given by $\exp(sM')\in\GL(V)$ via the
identifications $V \stackrel{\sim}{\rightarrow} M_{\gamma(0)}$ and
$V \stackrel{\sim}{\rightarrow} M_{\gamma(1)}$ chosen above. By $s$ we denote
here the unique real number such that $\tilde{\gamma}(t):= \tilde{x}_0
+ts\lambda \in \tilde{\Lambda}$ defines a lift of $\gamma$.

The shift functor on objects of $\mathcal{F}^0(E^\tau)$ is defined by 
$$(\Lambda,\alpha,M)[1] := (\Lambda,\alpha+1,M).$$

Let $(\Lambda_\nu,\alpha_\nu,M_\nu)$ with $\nu=1,2$ be two objects of
$\mathcal{F}^0(E^\tau)$. If $\Lambda_1\ne\Lambda_2$ we define 
\begin{multline}
 \Hom_{\mathcal{F}^0(E^\tau)}((\Lambda_1,\alpha_1,M_1),
 (\Lambda_2,\alpha_2,M_2)) = \notag\\ =
\begin{cases}
   0& \text{if $\alpha_2-\alpha_1\notin[0,1)$},\\
   \bigoplus_{x\in\Lambda_1\cap\Lambda_2} \Hom((M_1)_x,(M_2)_x)&
   \text{if $\alpha_2-\alpha_1\in[0,1)$}.
\end{cases}  
\end{multline}

If, however, $\Lambda_1=\Lambda_2$, we automatically have
$\alpha_1-\alpha_2\in\mathbb{Z}$ and we define
\begin{multline}
\Hom_{\mathcal{F}^0(E^\tau)}((\Lambda_1,\alpha_1,M_1),
(\Lambda_2,\alpha_2,M_2)) = \notag\\ =
\begin{cases}
  0& \text{if $\alpha_2-\alpha_1\notin\{0,1\}$},\\
  H^0(\Lambda_1,\HOM(M_1,M_2))& \text{if
  $\alpha_2=\alpha_1$},\\
  H^1(\Lambda_1,\HOM(M_1,M_2))& \text{if
  $\alpha_2=\alpha_1+1$}.
\end{cases}
\end{multline}

Here we denote by $\HOM(M_1,M_2)$ the sheaf of
homomorphisms of the two local systems defined by the vector spaces $V_\nu$ and
the automorphisms $M_\nu\in\GL(V_\nu)$. The vector space for this local system
is $V=\Hom(V_1,V_2)$ and the automorphism $M\in\GL(V)$ is given by
$M(f)=M_2\circ f\circ M_1^{-1}$ for $f\in V$, because we choose always the same
orientation on $\Lambda_1$.

Since $\Lambda_1$ is a circle, it is easy to compute the sheaf cohomology
appearing in the definition. The result is the following:
$H^0(\Lambda_1,M) \cong \ker(M-\mathbf{1}_V)$ and $H^1(\Lambda_1,M) \cong
\coker(M-\mathbf{1}_V)$. This gives, in particular, 
$$H^0(\Lambda_1, \HOM(M_1,M_2)) \cong \{f\in \Hom(V_1,V_2)\mid M_2\circ f=
f\circ M_1\}$$ and 
$$H^1(\Lambda_1, \HOM(M_1,M_2)) \cong
\Hom(V_1,V_2)/M_2\circ\Hom(V_1,V_2)\circ M_1^{-1}.$$

As an application  we obtain for $k=0,1$ 
$$H^k(\Lambda,\HOM(M_1,M_2)) \ne 0 \iff \text{$M_1$ and $M_2$ have a common
  eigenvalue}.$$

Since $\ker(M-\mathbf{1}_V)^\ast \cong \coker({}^tM-\mathbf{1}_V) =
\coker({}^tM^{-1}-\mathbf{1}_V)$ we obtain a canonical isomorphism
$H^0(\Lambda_1,M)^\ast \cong H^1(\Lambda_1, M^\vee)$ (where ${}^\ast$ denotes
the dual of a vector space, ${}^tM$ the transposed homomorphism and $M^\vee$
the dual local system, which is given by the automorphism ${}^tM^{-1}$.)

Combining this with our definitions we obtain the following ``Symplectic Serre
Duality'':

\begin{lemma}\label{symserre}
   If $(\Lambda_i,\alpha_i,M_i)$ are objects in $\mathcal{F}^0(E^\tau)$ then
   there exists a canonical isomorphism
   $$\Hom((\Lambda_1,\alpha_1,M_1),(\Lambda_2,\alpha_2,M_2)[1]) \cong
   \Hom((\Lambda_2,\alpha_2,M_2),(\Lambda_1,\alpha_1,M_1))^\ast.$$
\end{lemma}

To complete the definition of the category $\mathcal{F}^0(E^\tau)$ we have to
define the composition of morphisms. This will be the first place where
dependence on $\tau$ occurs. Because this is not yet done, we cannot speak
about functoriality in Lemma \ref{symserre}. 

Let $(\Lambda_1,\alpha_1,M_1) \stackrel{u}{\rightarrow}
(\Lambda_2,\alpha_2,M_2) \stackrel{v}{\rightarrow} (\Lambda_3,\alpha_3,M_3)$ be
two non-zero morphisms in $\mathcal{F}^0(E^\tau)$. 
The definitions imply $\alpha_1\le \alpha_2\le \alpha_3, \alpha_2\le
\alpha_1+1$ and $\alpha_3\le \alpha_2+1$. To define $v\circ u$ we have to
distinguish several cases. 

\textsl{Case (i):} $\alpha_1< \alpha_2 < \alpha_3 < \alpha_1+1$.

\nopagebreak
This is the
case with three different underlying Lagrangians. Only this case is considered
in \cite{PZ} and their definition is the following: 

Assume 
$u\in \Hom((M_1)_{x_1}, (M_2)_{x_1}) \subset
\Hom((\Lambda_1,\alpha_1,M_1),(\Lambda_2, \alpha_2,M_2))$ and $v\in
\Hom((M_2)_{x_2}, (M_3)_{x_2})$ with $x_1\in \Lambda_1\cap \Lambda_2, x_2\in
\Lambda_2\cap \Lambda_3$. For any $x_3\in \Lambda_1\cap \Lambda_3$, the
component of $v\circ u$ in $\Hom((M_1)_{x_3}, (M_3)_{x_3})$ will be given by
the following expression:
\begin{equation}
  \label{defmult}
  \sum_\phi \exp(2\pi i\int\phi^\ast(\tau dx\wedge dy))P(M_3)\circ v \circ
  P(M_2) \circ u \circ P(M_1)  
\end{equation}
where the sum is performed over  equivalence classes of quadruples
$(\phi; z_1, z_2, z_3)$ where $\phi:D\rightarrow E^\tau$ is a holomorphic map,
$D\subset \mathbb{C}$ denotes the unit disc and $\{z_1, z_2, z_3\}$ are
distinct points on the boundary of $D$ such that the order $(z_3, z_2, z_1)$
coincides with the usual orientation of $\mathbb{C}$. We require $\phi(z_\nu) =
x_\nu$ and that $\phi$ maps the arc connecting $z_\nu$ with $z_{\nu-1}$ to
$\Lambda_\nu$. Equivalent quadruples are obtained by applying automorphisms of
$D$ to $\phi$ and the $z_\nu$. By $P(M_\nu): (M_\nu)_{x_{\nu-1}} \rightarrow
(M_\nu)_{x_{\nu}}$ we denote the parallel transport on the local system
$M_\nu$. Using the identification $(M_\nu)_{x_\nu} \cong V_\nu$ chosen above
and assuming $M_\nu = \exp(M_\nu') \in \GL(V_\nu)$, the map $P(M_\nu)$ is
identified with $\exp(sM_\nu') \in \GL(V_\nu)$ where $s$ is the unique real
number such that $[0,1] \ni t \mapsto \tilde{x}_{\nu-1} + ts\lambda_\nu \in
\mathbb{R}^2$ is a lift to $\mathbb{R}^2$ of the path connecting $x_{\nu-1}$
with $x_\nu$ which is the image under $\phi$ of the corresponding arc on the
boundary of $D$. The weights in the sum can be computed as $\exp(2\pi i\tau
A_\phi)$ where $A_\phi$ denotes the Euclidean area of a triangle in
$\mathbb{R}^2$ whose image in $E^\tau$ is $\phi(D)$. (The edges of this
triangle can be used as the lifts of the paths connecting the $x_\nu$ as
above.) 

\textsl{Case (ii):} $\alpha_3>\alpha_1+1$.
 
We define $v\circ u=0$.

\textsl{Case (iii):} Precisely two of the $\alpha_\nu$ coincide and
$\alpha_1+1<\alpha_3$.

If $\alpha_1=\alpha_2<\alpha_3$ we have $\Lambda_1=\Lambda_2\ne \Lambda_3$. Let
us abbreviate $\boldsymbol{\Lambda}_\nu=(\Lambda_\nu,\alpha_\nu,M_\nu)$ and
define the composition via the commutative diagram

$$
\begin{CD}
  \Hom(\boldsymbol{\Lambda}_1,\boldsymbol{\Lambda}_2) \otimes
  \Hom(\boldsymbol{\Lambda}_2,\boldsymbol{\Lambda}_3) @>>>
  \Hom(\boldsymbol{\Lambda}_1,\boldsymbol{\Lambda}_3) \\
  @VV=V @VV=V\\
  H^0(\Lambda_1, \HOM(M_1,M_2))\otimes
  \bigoplus\limits_{x} \Hom((M_2)_x, (M_3)_x) @>>> 
  \bigoplus\limits_{x} \Hom((M_1)_x, (M_3)_x)\\ 
  \\ 
\end{CD}
$$
The sums in the bottom row are performed over all $x\in\Lambda_1\cap\Lambda_3$
and the morphism sends $\varphi\otimes (f_x)_x$ to $(f_x\circ\varphi_x)_x$
where $\varphi_x: (M_1)_x\rightarrow (M_2)_x$ denotes the map on stalks induced
by $\varphi$. The definition in the case $\alpha_1<\alpha_2=\alpha_3$ is
similar. 

\textsl{Case (iv):} $\alpha_1=\alpha_2=\alpha_3$. 

Here we have $\Lambda_1 =
\Lambda_2 =\Lambda_3$ and we can use the composition of homomorphisms of local
systems to define the composition in $\mathcal{F}^0(E^\tau)$.

\textsl{Case (v):} We are not in case (ii) and two of the $\alpha_\nu$ differ
by $1$ or equivalently $\alpha_3=\alpha_1+1$.

By using the canonical
isomorphism $\Hom(U\otimes V^\ast, W^\ast) = \Hom(W\otimes U, V)$ (where
$U,V,W$ are vector spaces) and symplectic Serre duality we reduce this
situation to other cases as follows: If $\alpha_1<\alpha_2<\alpha_3=\alpha_1+1$
we are led by the diagram
$$
\begin{CD}
  \Hom(\boldsymbol{\Lambda}_1,\boldsymbol{\Lambda}_2) \otimes 
  \Hom(\boldsymbol{\Lambda}_2,\boldsymbol{\Lambda}_3) @>>>
  \Hom(\boldsymbol{\Lambda}_1,\boldsymbol{\Lambda}_3) \\
  @VV=V @VV=V\\
  \Hom(\boldsymbol{\Lambda}_1,\boldsymbol{\Lambda}_2) \otimes
  \Hom(\boldsymbol{\Lambda}_3[-1],\boldsymbol{\Lambda}_2)^\ast @>>> 
  \Hom(\boldsymbol{\Lambda}_3[-1],\boldsymbol{\Lambda}_1)^\ast   \\
\end{CD}
$$
to a composition of the following kind
$$
\begin{CD}
  \Hom(\boldsymbol{\Lambda}_3[-1],\boldsymbol{\Lambda}_1) \otimes
  \Hom(\boldsymbol{\Lambda}_1,\boldsymbol{\Lambda}_2)  @>>> 
  \Hom(\boldsymbol{\Lambda}_3[-1],\boldsymbol{\Lambda}_2) \\
\end{CD}
$$
which was studied in case (iii). If $\alpha_2=\alpha_1+1$ we have
$\alpha_3=\alpha_2$ (otherwise we would be in case (ii)) and thus
$\Lambda_1=\Lambda_2 = \Lambda_3$. Similarly, if $\alpha_3 = \alpha_2 +1$ we
have $\alpha_1 = \alpha_2$ and $\Lambda_1=\Lambda_2 = \Lambda_3$. With the aid
of Lemma \ref{symserre} both cases are reduced to case (iv).

\begin{rem}
  Since symplectic Serre duality does not involve $\tau$, the definition of the
  composition of morphisms depends on $\tau$ in case (i) only.
\end{rem}

\begin{rem}
  The definition, especially case (v), is made precisely to get functoriality
  in Lemma \ref{symserre}.
\end{rem}

Because the derived category is additive and contains in particular finite
direct sums and a zero object, we would have no chance to find an equivalence
of categories if we don't have direct sums in Fukaya's category. By
construction, the category $\mathcal{F}^0(E^\tau)$ is an \textbf{Ab}-category
(also called a preadditive category), but it does not contain all direct sums
(biproducts). To see this let $(\Lambda_\nu, \alpha_\nu, M_\nu)$ be two objects
of $\mathcal{F}^0(E^\tau)$ and assume their direct sum $(\Lambda,\alpha,M) =
(\Lambda_1,\alpha_1,M_1)\oplus (\Lambda_2,\alpha_2,M_2)$ exists in
$\mathcal{F}^0(E^\tau)$. Then we have projections $p_\nu:(\Lambda,\alpha,M)
\rightarrow (\Lambda_\nu, \alpha_\nu, M_\nu)$ and embeddings
$i_\nu:(\Lambda_\nu, \alpha_\nu, M_\nu) \rightarrow (\Lambda,\alpha,M)$
fulfilling $p_\nu i_\nu = \mathbf{1}_{(\Lambda_\nu,\alpha_\nu,M_\nu)}$ and $i_1
p_1 + i_2 p_2 = \mathbf{1}_{(\Lambda,\alpha,M)}$. Assuming the
$(\Lambda_\nu,\alpha_\nu,M_\nu)$ are not zero objects, we obtain that neither
$p_\nu$ nor $i_\nu$ is the zero homomorphism. 

If we would have $\Lambda_1\ne \Lambda_2$, using 
$\Hom((\Lambda_1,\alpha_1,M_1),(\Lambda_2,\alpha_2,M_2)) \ne 0$
we would obtain directly from the definitions 
$\Hom((\Lambda_2,\alpha_2,M_2),(\Lambda_1,\alpha_1,M_1)) = 0$. Hence
$\Lambda_1 = \Lambda_2$ and we conclude that $\mathcal{F}^0(E^\tau)$
contains the direct sum of the objects as above only if $\Lambda_1=\Lambda_2$.

We use the following general construction to enlarge our category to
one with direct sums. 

Let $\mathcal{A}$ be an \textbf{Ab}-category. We define the category
$\underline{\mathcal{A}}$ to have as objects the ordered $k$-tuples of objects
of $\mathcal{A}$ where $k\ge 0$ runs through the integers. The morphisms are
formed by matrices of morphisms of $\mathcal{A}$ and their composition is given
by usual matrix multiplication. More formally $\ob(\underline{\mathcal{A}}) :=
\coprod\limits_{k\ge0} \prod\limits_{\nu=1}^k\ob(\mathcal{A})$. A $0$-tuple is
here considered to be a unique object of $\underline{\mathcal{A}}$ serving as a
zero object. We denote it by $\underline{0}$. If $\underline{A}=(A_1,\ldots,
A_k)$ and $\underline{B}=(B_1,\ldots, B_l)$ are objects of
$\underline{\mathcal{A}}$ with $k,l>0$, we define
$$\Hom_{\underline{\mathcal{A}}}(\underline{A}, \underline{B}) := \prod_{(i,j)}
\Hom_{\mathcal{A}}(A_i, B_j).$$ If $\underline{A}=\underline{0}$ or
$\underline{B}=\underline{0}$, we define
$\Hom_{\underline{\mathcal{A}}}(\underline{A}, \underline{B}) := 0$. This makes
$\underline{0}$ a zero object (this means it is initial and terminal). If the
category $\mathcal{A}$ contains already a zero object $0$, then all $k$-tuples
($k\ge 0$) of the form $(0,\ldots,0)$ are isomorphic to $\underline{0}$. The
composition of morphisms in $\underline{\mathcal{A}}$ is defined by the usual
formula for matrix multiplication. The composition $\underline{A} \rightarrow
\underline{0} \rightarrow \underline{B}$ is by definition the morphism whose
components are zero morphisms. It is easy to see that $\underline{\mathcal{A}}$
is an additive category and that the additive functor $\mathcal{A} \rightarrow
\underline{\mathcal{A}}$ sending an object to the $1$-tuple consisting of the
same object, is a universal functor from $\mathcal{A}$ to an additive category.

We apply this construction to $\mathcal{A}=\mathcal{F}^0$ and define the
(cohomology of the) \emph{Fukaya--Kontsevich} category
$\mathcal{FK}^0:=\underline{\mathcal{F}^0}$. This is an additive category and
the correctly formulated result from A. Polishchuk and E. Zaslow \cite{PZ} is
Theorem \ref{equiv}. 

Now, having the correct category, we are going to define some functors on
it. The study of $\mathcal{FK}^0(E^\tau)$ requires to consider the complex tori
$E^{r\tau}$ for all positive integers $r$. We denote by
$p_r:E^{r\tau}\rightarrow E^\tau$ the map defined by $p_r(x,y) = (rx,y)$. It
corresponds somehow to the isogeny $\pi_r:E_{r\tau}\rightarrow E_\tau$ on the
holomorphic side of the mirror story. Next we define additive functors
$p_{r\ast}:\mathcal{FK}^0(E^{r\tau})\rightarrow\mathcal{FK}^0(E^\tau)$ and 
$p_r^\ast:\mathcal{FK}^0(E^{\tau})\rightarrow\mathcal{FK}^0(E^{r\tau})$ which
are compatible with shifts. They are analogous to the functors $\pi_{r\ast}$
and $\pi_r^\ast$. 

Let $(\Lambda,\alpha,M)$ be an object in $\mathcal{FK}^0(E^{r\tau})$. We define
$$p_{r\ast}(\Lambda,\alpha,M) := (p_r(\Lambda), \alpha',p_{r\ast}M)$$ where
$\alpha'$ is the unique possible value lying in the same interval
$(k-\frac{1}{2},k+\frac{1}{2}]$ with $k\in\mathbb{Z}$ as $\alpha$ lies and
$p_{r\ast}M$ is the direct image of the local system $M$ in the sense of
sheaves of complex vector spaces. This means, if $p_r:\Lambda\rightarrow
p_r(\Lambda)$ is of degree $d$ and $M\in\GL(V)$ then
$p_{r\ast}M\in\GL(V^{\oplus d})$ is given by $(p_{r\ast}M)(v_1,\ldots,v_d) =
(v_2,\ldots,v_{d},Mv_1)$.

If $p_r(\Lambda_1)\ne p_r(\Lambda_2)$ we obtain 
$$p_{r\ast}:
\Hom((\Lambda_1,\alpha_1,M_1),(\Lambda_2,\alpha_2,M_2)) \rightarrow
\Hom(p_{r\ast}(\Lambda_1,\alpha_1,M_1),p_{r\ast}(\Lambda_2,\alpha_2,M_2))$$ 
in an obvious way from a bijection between the direct summands of the stalk
$(p_{r\ast}M_\nu)_x =
\bigoplus_{y\in\Lambda_\nu, p_r(y)=x}(M_\nu)_y$ with the stalks of
$M_\nu$ at the preimages of $x$ on $\Lambda_\nu$. 

If $\Lambda_1=\Lambda_2$, we use the canonical homomorphism of sheaves
$$p_{r\ast}\HOM(M_1,M_2) \rightarrow\HOM(p_{r\ast}M_1,p_{r\ast}M_2)$$ and the
fact that $p_r$ is a local homeomorphism to get the required map.

If $\Lambda_1\ne\Lambda_2$ but $p_r(\Lambda_1)=p_r(\Lambda_2)$ we necessarily
have $\Lambda_1\cap\Lambda_2=\emptyset$, hence $p_{r\ast}$ is zero in this
case. 

To verify compatibility of the functor $p_{r\ast}$ with compositions, we
consider three 
objects $(\Lambda_\nu,\alpha_\nu,M_\nu)$ in $\mathcal{FK}^0(E^{r\tau})$. If at
least two of the $\Lambda_\nu$ coincide, compatibility is easily obtained from
the definitions. If $\Lambda_1\ne \Lambda_2\ne \Lambda_3 \ne \Lambda_1$ we have
to compare two sums over certain holomorphic maps $\phi^\tau:D\rightarrow
E^\tau$ and $\phi^{r\tau}:D\rightarrow E^{r\tau}$ respectively (see
(\ref{defmult})). The map $p_r$ defines a bijection between the images of the
lifts of these maps to $\mathbb{R}^2$. These images are triangles and their
Euclidean areas $A_{\phi^\tau}$ and $A_{\phi^{r\tau}}$ fulfill $A_{\phi^\tau} =
rA_{\phi^{r\tau}}$. This implies the equality we want.

Let now $(\Lambda,\alpha,M)$ be an object in $\mathcal{FK}^0(E^{\tau})$. Assume
the preimage $p_r^{-1}(\Lambda)$ consists of $n$ connected components
$\Lambda^{(1)},\ldots,\Lambda^{(n)}$. Then the restrictions
$p_r^{(k)}:\Lambda^{(k)}\rightarrow\Lambda$ are of degree $d:=r/n$ for all
$k$. We define
$$p_r^\ast(\Lambda,\alpha,M) := \bigoplus\limits_{k=1}^n
(\Lambda^{(k)},\alpha',(p_r^{(k)})^{\ast}M),$$ where $\alpha'$ is in the same
interval $(k-\frac{1}{2},k+\frac{1}{2}]$ with $k\in\mathbb{Z}$ as $\alpha$ is,
and $(p_r^{(k)})^{\ast}M$ denotes the sheaf on $\Lambda^{(k)}$ obtained by
pulling back the local system $M$ considered as a locally constant sheaf of
complex vector spaces. This means more explicitly that this local system is
defined by the vector space $V$ and the automorphism $M^d\in\GL(V)$.  Observe
that $p_r^\ast$ does not have values in $\mathcal{F}^0(E^\tau)$, because we
need direct sums if $n>1$.

To define $p_r^\ast$ on morphisms we have to distinguish two
cases. In the first case we assume $\Lambda:=\Lambda_1=\Lambda_2$. From sheaf
theory we know that there is a canonical homomorphism
$\HOM((p_r^{(k)})^{\ast}M_1, (p_r^{(k)})^{\ast}M_2) \rightarrow
(p_r^{(k)})^{\ast}\HOM(M_1,M_2)$. By taking cohomology, we obtain $p_r^\ast$.

In the second case we assume $\Lambda_1\ne\Lambda_2$. The definition is now
straightforward: If $f \in
\Hom((\Lambda_1,\alpha_1,M_1),(\Lambda_2,\alpha_2,M_2)) =
\bigoplus\limits_{x\in\Lambda_1\cap\Lambda_2} \Hom((M_1)_x,(M_2)_x)$ has
components $f_x \in \Hom((M_1)_x,(M_2)_x)$ with ${x\in\Lambda_1\cap\Lambda_2}$,
we define the component $(p_r^\ast f)_y$ of $p_r^\ast f$ corresponding to $y\in
p_r^{-1}\Lambda_1\cap p_r^{-1}\Lambda_2$ by $$(p_r^\ast f)_y = f_{p_r(y)}\in
\Hom((M_1)_{p_r(y)},(M_2)_{p_r(y)}).$$ We used here the canonical isomorphism
$((p_r^{(k)})^{\ast}M)_y = M_{p_r(y)}$. 

A similar reasoning as for the functor $p_{r\ast}$ leads to compatibility of
$p_r^\ast$ with compositions. This establishes that we really defined two
functors $p_{r\ast}$ and $p_r^\ast$.

We shall need also the pull-back under translations. If $(x_0,y_0) \in
\mathbb{R}^2$ we let $t(x,y) = (x-x_0, y-y_0)$ be a translation on $E^\tau$. We
define $t^\ast(\Lambda,\alpha,M) := (t^{-1}(\Lambda),\alpha,t^\ast M)$. Since
$t$ is an 
isomorphism the corresponding map $t^\ast$ on morphisms and compatibility with
composition is obvious. Similarly, we define $t_\ast$. Directly from the
definitions we obtain $(t\circ p_r)_\ast = t_\ast\circ p_{r\ast}$ and $(t\circ
p_r)^\ast = p_r^\ast\circ t^\ast$.

\begin{lemma}\label{sympladj}
  Let $p_r: E^{r\tau}\rightarrow E^\tau$ be as above and $t:E^\tau \rightarrow
  E^\tau$ a translation of the form $t(x,y) = (x+\frac{m}{n},y)$, with $m,n$
  integers. Define $p=t\circ p_r: E^{r\tau} \rightarrow E^\tau$. Let
  $(\Lambda_1,\alpha_1,M_1), (\Lambda_2,\alpha_2,M_2)$ be objects in
  $\mathcal{FK}^0(E^\tau)$ and $\mathcal{FK}^0(E^{r\tau})$ respectively. Then
  there are functorial isomorphisms:
  $$\Hom(p^\ast(\Lambda_1,\alpha_1,M_1),(\Lambda_2,\alpha_2,M_2)) \cong
  \Hom((\Lambda_1,\alpha_1,M_1),p_{\ast}(\Lambda_2,\alpha_2,M_2))$$
  and
  $$\Hom(p_{\ast}(\Lambda_2,\alpha_2,M_2),(\Lambda_1,\alpha_1,M_1)) \cong
   \Hom((\Lambda_2,\alpha_2,M_2),p^\ast(\Lambda_1,\alpha_1,M_1)).$$
\end{lemma}

\begin{proof}
Since $(t\circ p_r)^\ast = p_r^\ast \circ t^\ast$ and $(t\circ
p_r)_\ast = t_\ast \circ p_{r\ast}$ it is enough to prove the claim in case
$p=p_r$ and $p=t$ separately. But in the case $p=t$ the statements are
obvious. For the rest of the proof we assume $p=p_r$. 

If $\Lambda_1=p_r(\Lambda_2)$ the first statement is the usual
adjointness of $p^{\ast}$ and $p_\ast$ for sheaves of vector spaces. To obtain
the second statement in this case, we have to use that $p_r$ is a local
homeomorphism. 

If $\Lambda_1\ne p_r(\Lambda_2)$ we have
\begin{multline}
\Hom((\Lambda_1,\alpha_1,M_1),p_{r\ast}(\Lambda_2,\alpha_2,M_2)) = 
\bigoplus_{x\in\Lambda_1\cap p_r(\Lambda_2)} \Hom((M_1)_x, p_{r\ast}(M_2)_x)
=\notag \\
\bigoplus_{\substack{x\in\Lambda_1\cap p_r(\Lambda_2)\\y\in\Lambda_2,
    p_r(y)=x}} 
  \Hom((M_1)_x, (M_2)_y) =  
\bigoplus_{y\in p_r^{-1}(\Lambda_1)\cap \Lambda_1} 
  \Hom((M_1)_{p_r(y)}, (M_2)_y) = \\
\bigoplus_{y\in p_r^{-1}(\Lambda_1)\cap \Lambda_1}
  \Hom((p_r^\ast M_1)_y, (M_2)_y) = 
\Hom(p_r^\ast(\Lambda_1,\alpha_1,M_1),(\Lambda_2,\alpha_2,M_2))
\end{multline}
where we used
the canonical isomorphisms $(p_{r\ast}M_2)_x = \bigoplus_{y\in\Lambda_2,
  p_r(y)=x} (M_2)_y$ and $(p_r^\ast M_1)_y = (M_1)_{p_r(y)}$. 

Suppose $(\Lambda_1', \alpha_1', M_1')$ is another object in
$\mathcal{FK}^0(E^\tau)$ with $\Lambda_1'\ne \Lambda_1$ and $\Lambda_1'\ne
p_r(\Lambda_2)$. Let a morphism $(\Lambda_1',\alpha_1',M_1') \rightarrow
(\Lambda_1,\alpha_1,M_1)$ be given and consider the induced diagram:

$$
\begin{CD}
  \Hom(p_r^\ast(\Lambda_1,\alpha_1,M_1),(\Lambda_2,\alpha_2,M_2)) @>>\cong>
  \Hom((\Lambda_1,\alpha_1,M_1),p_{r\ast}(\Lambda_2,\alpha_2,M_2))\\
  @VVV @VVV\\
  \Hom(p_r^\ast(\Lambda_1',\alpha_1',M_1'),(\Lambda_2,\alpha_2,M_2)) @>>\cong>
  \Hom((\Lambda_1',\alpha_1',M_1'),p_{r\ast}(\Lambda_2,\alpha_2,M_2)).
\end{CD}
$$

To get functoriality we have to show commutativity of this diagram. This
involves the definition of composition in $\mathcal{FK}^0$. Assume the given
morphism is represented simply by $f:(M_1')_x \rightarrow (M_1)_x$ where $x\in
\Lambda_1\cap \Lambda_1'$. The left vertical arrow sends the morphism $g:
(p_r^\ast M_1)_y = (M_1)_{p_r(y)} \rightarrow (M_2)_y$ with $y\in
p_r^{-1}(\Lambda_1)\cap \Lambda_2$ to $g\circ p_r^\ast(f) \in \bigoplus_{z\in
  \Lambda_2, p_r(z)=x} \Hom((p_r^\ast M_1')_z, (M_2)_z)$ whose component
corresponding to $z$ is $\sum_\phi \exp(2\pi ir\tau A_\phi) P(M_2)\circ g \circ
P(p_r^\ast M_1) \circ (p_r^\ast f)_z \circ P(p_r^\ast M_1')$. The sum is
performed over equivalence classes of certain maps $\phi:D\rightarrow
E^{r\tau}$ sending a boundary point to $z$ and $(p_r^\ast f)_z$ is obtained
from $f$ by composition with the identification $(p_r^\ast M_1')_z =
(M_1')_x$. 

Our isomorphism sends this to the morphism on the right hand side whose
component in 
\begin{multline}
  \Hom((\Lambda_1',\alpha_1',M_1'), p_r^\ast(\Lambda_2,\alpha_2,M_2)) =
  \notag\\ = \bigoplus_{x\in p_r(\Lambda_2)\cap \Lambda_1'} \Hom((M_1')_x,
  (p_{r\ast} M_2)_x) =\\= \bigoplus_{y\in\Lambda_2 \cap 
  p_r^{-1}(\Lambda_1')} \Hom((M_1')_{p_r(y)}, (M_2)_y)
\end{multline}

corresponding to $z\in
p_r^{-1}\Lambda_1'\cap \Lambda_2$ is given by the above sum. This coincides
with the sum we would write down for the composition on the right hand side,
because $p_r$ defines a bijection between the images of the lifts of the
holomorphic maps $\phi$ to the universal covers $\mathbb{R}^2$ of $E^\tau$ and
$E^{r\tau}$ respectively. These images are triangles and if we denote their
Euclidean area by $A_\phi^\tau$ and $A_\phi^{r\tau}$ we obtain from the
definition of $p_r$ the equation $A_\phi^\tau = rA_\phi^{r\tau}$ ant this shows
that both sums have the same weights. 

To get functoriality in general we have to consider situations involving both
cases considered above. In all such situations the commutativity of the
corresponding diagram follows immediately from the definition of composition in
$\mathcal{FK}^0$. A similar consideration gives the second statement.
\end{proof}

\begin{rem}
  The same statement is true for $\pi_{r\ast}$ and $\pi_r^\ast$.
\end{rem}

To prepare the proof of Theorem \ref{equiv} let us consider the following
Cartesian diagram 
$$
\begin{CD}
  \tilde{E} @>>{\tilde{p}_{r_1}}> E^{r_2\tau} \\
  @V{\tilde{p}_{r_2}}VV @VV{p_{r_2}}V\\
  E^{r_1\tau} @>>{p_{r_1}}> E^\tau
\end{CD}
$$
where $\tilde{E} = E^{r\tau}\times \mathbb{Z}/d\mathbb{Z}$ with
$d=\gcd(r_1,r_2)$ and $r=\frac{r_1 r_2}{d}$. We denote by $\tilde{p}_{r_i,\nu}$
the restriction of $\tilde{p}_{r_i}$ to $E^{r\tau}\times\{\nu\}$, which is the
composition of $p_{\frac{r_{3-i}}{d}}: E^{r\tau} \rightarrow E^{r_i\tau}$ with
a translation of the form $(x,y)\mapsto (x-n,y)$ on $E^{r_i\tau}$. In our
application these translations are determined by corresponding translations on
the elliptic curves $E_{r_i\tau}$. Under these assumptions the following lemma
holds: 

\begin{lemma}\label{symplbasechange}
  There exists a functorial isomorphism
  \begin{multline}
  \Hom(p_{r_1\ast}(\Lambda_1,\alpha_1,M_1),
  p_{r_2\ast}(\Lambda_2,\alpha_2,M_2)) \cong\notag \\ \cong \bigoplus
  \Hom(\tilde{p}_{r_2,\nu}^\ast(\Lambda_1,\alpha_1,M_1),
  \tilde{p}_{r_1,\nu}^\ast(\Lambda_2,\alpha_2,M_2))    
  \end{multline}
\end{lemma}
\begin{proof}
If we apply Lemma \ref{sympladj} to $p_{r_1}$ we obtain a
functorial isomorphism: 
\begin{multline}
  \Hom(p_{r_1\ast}(\Lambda_1,\alpha_1,M_1),
  p_{r_2\ast}(\Lambda_2,\alpha_2,M_2)) \cong\notag \\ \cong 
  \Hom((\Lambda_1,\alpha_1,M_1),
  p_{r_1}^\ast p_{r_2\ast}(\Lambda_2,\alpha_2,M_2)).   
  \end{multline}
On the other hand, by applying Lemma \ref{sympladj} to the maps
$\tilde{p}_{r_{1,\nu}}$ we obtain:
\begin{multline}
  \Hom(\tilde{p}_{r_2,\nu}^\ast(\Lambda_1,\alpha_1,M_1),
  \tilde{p}_{r_1,\nu}^\ast(\Lambda_2,\alpha_2,M_2))
  \cong\notag \\ \cong
  \bigoplus_{\nu} \Hom((\Lambda_1,\alpha_1,M_1),
  \tilde{p}_{r_2,\nu\ast} \tilde{p}_{r_1,\nu}^\ast(\Lambda_2,\alpha_2,M_2))
  \cong \\ \cong
  \Hom((\Lambda_1,\alpha_1,M_1),
  \tilde{p}_{r_2\ast}\tilde{p}_{r_1}^\ast (\Lambda_2,\alpha_2,M_2)).   
\end{multline}
To get the claim we need to see that the functors $\tilde{p}_{r_2\ast} \circ
\tilde{p}_{r_1}^\ast$ and $p_{r_1}^\ast\circ p_{r_2\ast}$ are isomorphic, but
this follows easily from the definitions because $p_{r_1}$ and
$\tilde{p}_{r_1}$ are local homeomorphisms.
\end{proof}

\section{The equivalence}
\label{sec:equiv}

\begin{thm}\label{equiv}
  $\Phi_\tau: \mathcal{D}^b(E_\tau) \rightarrow \mathcal{FK}^0(E^\tau)$ is an
  equivalence of additive categories compatible with the shift functors. 
\end{thm}

The rest of this section consists of the proof of this theorem.

Recall (Proposition \ref{subcat}) that we replace $\mathcal{D}^b(E_\tau)$ by
the equivalent subcategory whose objects are finite direct sums of shifted
indecomposable coherent sheaves on $E_\tau$.

To define $\Phi=\Phi_\tau$ on objects it suffices to define the objects
$\Phi(A)$ where 
$A=A[0]$ is an indecomposable coherent sheaf on $E_\tau$. 

We repeat here the definition given in \cite{PZ}. The strategy is the
following. First they define $\Phi_{r\tau}(A)$ for all positive integers $r$
and all vector bundles $A$ being the tensor product of a line bundle with a
unipotent bundle on $E_\tau$. In the second step they extend this to all
indecomposable vector bundles by using Theorem \ref{vb}. Finally, they write
down the definition of $\Phi_\tau(A)$ for torsion sheaves $A$. 
These definitions are more or less dictated by the description of mirror
symmetry as T-duality \cite{SYZ}, \cite{Fuk}.
The definitions are the following: 

If $A=\mathcal{L}(\varphi)\otimes F(V,\exp(N))$ with
$\varphi=t^\ast_{a\tau+b}\varphi_0^{\phantom{1}}\cdot\varphi_0^{n-1}$  and $V$
a finite dimensional complex vector space and $N\in\End(V)$ a cyclic nilpotent
endomorphism (we call a nilpotent endomorphism cyclic, if the corresponding
$\mathbb{C}[T]$-module structure on $V$ is cyclic, which is equivalent to
$\dim\ker N=1$), we define
$$\Phi(A) = (\Lambda, \alpha, M)$$
where $\Lambda$ is given by $(a+t, (n-1)a+nt)$, $\alpha$ is the unique real
number with $-\frac{1}{2}<\alpha\le \frac{1}{2}$ and $\exp(i\pi\alpha) =
\frac{1+in}{\sqrt{1+n^2}}$ (i.e. $i\pi\alpha$ is a logarithm of the slope of
$\Lambda$) and $M= \exp(-2\pi ib\mathbf{1}_V + N)$.

As a first and useful observation we note here that 
$${t'}^{\ast}\Phi(A) = \Phi(t^\ast A)$$
if $A$ is a vector bundle on $E_\tau$ as before, $t:E_\tau\rightarrow E_\tau$
denotes the translation by $\frac{m}{n}\tau$ (i.e. $t(z)=z+\frac{m}{n}\tau$)
with integers $m,n$ and $t':E^\tau\rightarrow E^\tau$ is the corresponding
translation on the symplectic side, given by $t'(x,y) = (x-\frac{m}{n},y)$.

If $\pi_r: E_{r\tau}\rightarrow E_\tau$ is an isogeny and $A$ is a vector
bundle on $E_{r\tau}$ of the form considered above and
$p_r:E^{r\tau}\rightarrow E^\tau$ is the corresponding morphism on the
symplectic side, we define 
$$\Phi_\tau(\pi_{r\ast}A) := p_{r\ast}\Phi_{r\tau}(A).$$

More explicitly, this means for $A=\mathcal{L}(\varphi)\otimes F(V,\exp(N))$ as
above $\Phi_\tau(\pi_{r\ast}A) = (\Lambda, \alpha, M)$ where $\Lambda$ is given
by the line with $x$-intercept $\frac{ra}{n}$ and slope $\frac{n}{r}$, 
$-\frac{1}{2}<\alpha\le\frac{1}{2}$ is the corresponding logarithm and 
$M=\exp(-2\pi ib\mathbf{1}_V+N)$ because $p_r$ defines an isomorphism between
the underlying Lagrangian submanifolds.

Finally, let $A=S(a\tau+b, V, N)$ be a torsion sheaf on $E_\tau$ supported at
$x=a\tau+b$. Then $\Phi_\tau(A) := (\Lambda, \frac{1}{2}, 
\exp(2\pi ib\mathbf{1}_V+N))$ where $\Lambda$ is defined by $(-a, t)$.

To get a functor $\Phi$ we have to define its effect on morphisms. Using
compatibility with finite direct sums and with shifts, this amounts to define
$$\Phi_\tau: \Hom_{\mathcal{D}^b(E_\tau)} (A_1,A_2[n]) \rightarrow
\Hom_{\mathcal{FK}^0(E^\tau)} (\Phi_\tau(A_1), \Phi_\tau(A_2)[n]).$$ Since both
sides vanish if $n\notin \{0,1\}$, we are concerned with $n=0$ and $n=1$ only.
Using Serre duality (Lemmas \ref{ellserre} and \ref{symserre}) and
$\Hom_{\mathcal{D}^b(E_\tau)}(A_1, A_2[1]) = \Ext^1(A_1, A_2)$ we obtain the
isomorphism $\Phi_\tau$ for all $n$ from its definition in case $n=0$.

In the paper \cite{PZ} the considerations are restricted to those cases where
$\Phi(A_1)$ and $\Phi(A_2)$ have distinct underlying Lagrangian
submanifolds. We shall extend the definition to all cases. 

Let us repeat the definition of \cite{PZ} in the situation $A_i =
\mathcal{L}(\varphi_i)\otimes F(V_i, N_i)$ and $\Phi_\tau(A_i) = (\Lambda_i,
\alpha_i, M_i)$ with $\Lambda_1\ne \Lambda_2$. As usual, we write $\varphi_i =
t^\ast_{x_i}\varphi_0^{\phantom{1}}\cdot\varphi_0^{n_i-1}$ with
$x_i=a_i\tau+b_i$. Using the isomorphism  
$H^0(\mathcal{L}(\varphi))\otimes V\rightarrow H^0(\mathcal{L}(\varphi)\otimes
F(V, \exp(N)))$ explicitly given in \cite[Proposition 2]{PZ} and assuming
$n_1 < n_2$, we obtain an isomorphism  
$\Hom(A_1, A_2) \cong H^0(\mathcal{L}(\varphi_2^{\phantom{1}}\varphi_1^{-1}))
\otimes \Hom(V_1, V_2)$. Because $\Phi$ becomes complex linear it suffices to
define the images of elements of the form $\theta\otimes f$ with
$f\in\Hom(V_1,V_2)$ and $\theta$ a basis element of
$H^0(\mathcal{L}(\varphi_2^{\phantom{1}}\varphi_1^{-1}))$. The basis of this
vector space we are using is formed by the theta functions $t^\ast_{x_{12}}
\theta\big[\frac{k}{n_2-n_1}, 0\big]((n_2-n_1)\tau, (n_2-n_1)z)$. Here we use
the abbreviation $x_{12}=\frac{x_2-x_1}{n_2-n_1}$ and $0\le k < n_2-n_1$. Their
number $(n_2-n_1)$ is equal to the number of intersection points of $\Lambda_1$
and $\Lambda_2$. The definition of the isomorphism 
$$\Phi_\tau: \Hom_{\mathcal{D}^b(E_\tau)} (A_1,A_2) \rightarrow
\Hom_{\mathcal{FK}^0(E^\tau)} (\Phi_\tau(A_1), \Phi_\tau(A_2))$$
is now obtained by
specifying a bijection between these theta functions and the intersection
points. The choice of \cite{PZ} is to send $t^\ast_{x_{12}}
\theta\big[\frac{k}{n_2-n_1}, 0\big]((n_2-n_1)\tau, (n_2-n_1)z)$ to the
point $e_k= (\frac{a_2-a_1+k}{n_2-n_1}, \frac{n_1a_2-n_2a_1+n_1k}{n_2-n_1})$.

If $n_1 \ge n_2$ but
$\mathcal{L}(\varphi_1) \not\cong \mathcal{L}(\varphi_2)$ both spaces of
homomorphisms are zero.

Let us extend the definition to the case $\Lambda_1=\Lambda_2$. Under this
assumption we have $a_1=a_2$ and $\deg\mathcal{L}(\varphi_1) =
\deg\mathcal{L}(\varphi_2)$. In particular,
$\mathcal{L}(\varphi_2^{\phantom{1}}\varphi_1^{-1})$ is of degree zero. Because
the structure sheaf is the only line bundle of degree zero with nontrivial 
sections, we have $H^0(\mathcal{L}(\varphi_2^{\phantom{1}}\varphi_1^{-1}))=0$
if $\varphi_1\ne \varphi_2$. But $\varphi_1 = \varphi_2$ is equivalent to $b_1
= b_2$ which in turn is equivalent to the existence of a common eigenvalue of
$M_1$ and $M_2$. This is equivalent to the non-vanishing of
$\Hom_{\mathcal{F}^0}((\Lambda_1,\alpha_1,M_1),(\Lambda_1,\alpha_1,M_2))$. So
the spaces of morphisms are zero on both sides, if $\varphi_1\ne \varphi_2$. 

If $\varphi_1 = \varphi_2$ the operators $M_i$ have the same eigenvalue
$\exp(-2\pi ib_1)$. This implies 
$H^0(\Lambda_1, \HOM(M_1, M_2)) = \{f:V_1\rightarrow V_2\mid M_1\circ f=f\circ
M_2\} = \{f:V_1\rightarrow V_2\mid \exp(N_1)\circ f=f\circ \exp(N_2)\}$. The
composition of all these identifications with Lemma \ref{flatsect} gives the
definition of $\Phi$.

The next step is the extension of the definition of $\Phi$ to homomorphisms
between arbitrary vector bundles on $E_\tau$. Let $r_1,r_2$ be two positive
integers (possibly equal to one) and
$\mathcal{E}_i=\mathcal{L}(\varphi_i)\otimes F(V_i,\exp(N_i))$ vector
bundles on $E_{r_i\tau}$. Define $E:= E_{r_1\tau}\times_{E_\tau}E_{r_2\tau}$
and denote the projections by $\tilde{\pi}_{r_i}:E\rightarrow E_{r_i\tau}$. The
curve $E$ is a disjoint union of $d=\gcd(r_1,r_2)$ elliptic curves and the
restrictions $\tilde{\pi}_{r_i,\nu}$ of the projections $\tilde{\pi}_{r_i}$
onto the components are isogenies composed with translations. More precisely,
$E = E_{r\tau}\times \mathbb{Z}/d\mathbb{Z}$ with $r:=\frac{r_1r_2}{d}$ and
$\tilde{\pi}_{r_i,\nu}: E_{r\tau} \times \{\nu\} \rightarrow E_{r_i\nu}$ is the
composition of the isogeny $E_{r\tau} \rightarrow E_{r_i\tau}$ defined by
$\Gamma_{r_i\tau} \subset \Gamma_{r\tau}$ with the translation by $\nu\tau$ on
$E_{r_1\tau}$ and with the identity on $E_{r_2\tau}$. We could use the
translation by $s_i\tau$ on $E_{r_i\tau}$ for any pair of integers $(s_1,s_2)$
with $s_1-s_2 \equiv\nu \mod d$. Two such choices differ by a translation on
$E_{r\tau}$ by an integer multiple of $\tau$.

Using the flat base change theorem and the adjointness properties of $\pi_\ast$
and $\pi^\ast$, we obtain a canonical isomorphism
$$\Hom(\pi_{r_1\ast}\mathcal{E}_1, \pi_{r_2\ast}\mathcal{E}_2) \cong
\Hom(\tilde{\pi}_{r_1}^\ast\mathcal{E}_1, \tilde{\pi}_{r_2}^\ast\mathcal{E}_2)
\cong \bigoplus\limits_{\nu=1}^d \Hom(\tilde{\pi}_{r_1,\nu}^\ast\mathcal{E}_1,
\tilde{\pi}_{r_2,\nu}^\ast\mathcal{E}_2).$$ Because the situation on the
symplectic side is similar (Lemma \ref{symplbasechange}), we can define
$\Phi_\tau$ by the following commutative diagram:

$$
\begin{CD}
  \Hom(\pi_{r_1\ast}\mathcal{E}_1, \pi_{r_2\ast}\mathcal{E}_2) @>\sim>>
  \bigoplus\Hom(\tilde{\pi}_{r_1,\nu}^\ast\mathcal{E}_1,
  \tilde{\pi}_{r_2,\nu}^\ast\mathcal{E}_2) \\
  @VV{\Phi_\tau}V  @VV{\oplus\Phi}V\\
  \Hom(\Phi(\pi_{r_1\ast}\mathcal{E}_1), \Phi(\pi_{r_2\ast}\mathcal{E}_2))
  && 
  \bigoplus\Hom(\Phi(\tilde{\pi}_{r_1,\nu}^\ast\mathcal{E}_1),
  \Phi(\tilde{\pi}_{r_2,\nu}^\ast\mathcal{E}_2)) \\
  @VV{\cong}V  @VV{\cong}V\\
  \Hom(p_{r_1\ast}\Phi(\mathcal{E}_1), p_{r_2\ast}\Phi(\mathcal{E}_2))
  @>\sim>> 
  \bigoplus\Hom(\tilde{p}_{r_1,\nu}^\ast\Phi(\mathcal{E}_1),
  \tilde{p}_{r_2,\nu}^\ast\Phi(\mathcal{E}_2)) \\
\end{CD}
$$

In addition to the definition of $\Phi(\mathcal{E}_i)$ we use
$\Phi_{r\tau}(\pi^\ast(\mathcal{E})) \cong p^\ast(\Phi_\tau(\mathcal{E}))$ if
$\mathcal{E}\cong \mathcal{L}(\varphi)\otimes F(V,\exp(N))$ (see
\cite[Prop. 4]{PZ}) for isogenies $\pi$ and the compatibility of translations
with $\Phi$ observed earlier.

We still have to deal with the cases where $A_1$ or $A_2$ is a torsion
sheaf. The easy observation is that $\Hom(A_1, A_2) = \Ext^1(A_2, A_1) = 0$ if
$A_1$ is a torsion sheaf and $A_2$ is locally free. This corresponds nicely
to our definition on the symplectic side where
$\Hom((\Lambda_1,\alpha_1,M_1),(\Lambda_2, \alpha_2, M_2))=0$ if $\alpha_1 >
\alpha_2$. So we have to investigate just the cases where $A_2 = S(a_2\tau+b_2,
V_2, N_2)$ is a torsion sheaf. If $A_1=\mathcal{L}(\varphi_1)\otimes F(V_1,
N_1)$ the definition of \cite{PZ} is the following:

There are canonical isomorphisms $\Hom(A_1,A_2) = \Hom(V_1,V_2) =
V_1^\ast\otimes V_2$ and $\Hom(\Phi(A_1),\Phi(A_2)) = \Hom(V_1, V_2) = V_1^\ast
\otimes V_2$ because $\Lambda_1\cap \Lambda_2$ is one point. The isomorphism
$\Phi_\tau: V_1^\ast \otimes V_2 \rightarrow V_1^\ast \otimes V_2$ is by
definition 
\begin{multline}
  \exp(-\pi i\tau(na_2^2 + 2a_1a_2) - 2\pi i(a_2b_1+a_1b_2+na_2b_2))
  \cdot\notag \\ \cdot
  \exp(-(a_1+na_2)\mathbf{1}_{V_1^\ast}\otimes N_2 + a_2 {}^tN_1\otimes
  \mathbf{1}_{V_2}). 
\end{multline}

To extend the definition to arbitrary indecomposable locally free sheaves
$A_1$, assume $A_1=\pi_{r\ast}\mathcal{E}_1$ with
$\mathcal{E}_1=\mathcal{L}(\varphi_1)\otimes F(V_1,\exp(N_1))$ and $A_2$ a
torsion sheaf as above. We define $\Phi_\tau$ by the commutative diagram:

$$
\begin{CD}
  \Hom(\mathcal{E}_1, \pi_r^\ast A_2) @>{\sim}>> \Hom(\pi_{r\ast}\mathcal{E}_1,
  A_2)\\
  @V{\cong}V{\Phi_{r\tau}}V @VV{\Phi_\tau}V\\
  \Hom(\Phi(\mathcal{E}_1), \Phi(\pi_r^\ast A_2)) @>{\sim}>>
  \Hom(\Phi(\pi_{r\ast}\mathcal{E}_1), \Phi(A_2))\\
  @V{\cong}VV @V{\cong}VV\\
  \Hom(\Phi(\mathcal{E}_1), p_r^\ast\Phi(A_2)) @>{\sim}>>
  \Hom(p_{r\ast}\Phi(\mathcal{E}_1), \Phi(A_2)).\\ 
\end{CD}
$$

Here we used the definition of $\Phi$ on objects, Lemma \ref{sympladj} and
$\Phi(\pi_r^\ast A_2) = p_r^\ast\Phi(A_2)$ which follows easily from the
definitions.

Note that we need here the existence
of direct sums in the category $\mathcal{FK}^0(E^{r\tau})$ because $\pi_r^\ast
A_2$ is not indecomposable, it is the direct sum of skyscraper sheaves
supported at all preimages of the support on $A_2$. 

Finally, we have to deal with two indecomposable torsion sheaves $A_\nu$. If
they have distinct support, then $\Hom(A_1, A_2) = \Ext^1(A_1, A_2) = 0$. With
$A_\nu = S(a_\nu\tau+b_\nu, V_\nu, N_\nu)$ we have $\Phi(A_\nu) = (\Lambda_\nu,
\frac{1}{2}, M_\nu)$ where $\Lambda_\nu$ is given by $(-a_\nu, t)$ and
$M_\nu=\exp(2\pi ib_\nu\mathbf{1}_{V_\nu}+N_\nu)$. If $a_1\ne a_2$ we get  
$\Lambda_1\cap\Lambda_2=\emptyset$ and $\Hom(\Phi(A_1), \Phi(A_2)[n]) =0$. If
$a_1=a_2$ but $b_1\ne b_2$ then $\Lambda_1 = \Lambda_2$ but $M_1$ and $M_2$
have no common eigenvalues, hence $H^n(\Lambda_1, \HOM(M_1, M_2)) = 0$ for
$n=0,1$. If finally $A_1$ and $A_2$ have the same support, we obtain $\Hom(A_1,
A_2) = \Hom_{\mathcal{O}_{E_\tau,x_1}}((V_1,N_1),(V_2,N_2)) =
\{f\in\Hom(V_1,V_2)\mid f\circ N_1 = N_2\circ F\} = \{ f\in\Hom(V_1,V_2)\mid
f\circ M_1 = M_2\circ f\} = H^0(\Lambda_1, \HOM(M_1,M_2)).$

This completes the definition of the functor $\Phi$. The most important part of
the proof is to show that $\Phi$ is indeed a functor. This amounts to showing
compatibility of $\Phi$ with compositions. 

We have to show the commutativity of the diagram 
$$
\begin{CD}
  \Hom(A_1,A_2[k]) \otimes \Hom(A_2[k],A_3[l]) @>>> \Hom(A_1,A_3[l])\\
  @VV{\Phi_\tau\otimes\Phi_\tau}V @VV{\Phi_\tau}V\\
  \Hom(\Phi(A_1),\Phi(A_2)[k])\otimes \Hom(\Phi(A_2)[k], \Phi(A_3)[l]) @>>>
  \Hom(\Phi(A_1), \Phi(A_3)[l])
\end{CD}
$$
where the horizontal arrows are the compositions in the respective categories
and $0\le k\le l\le 1$. We write $\Phi(A_\nu) =
(\Lambda_\nu,\alpha_\nu,M_\nu)$. If $l=1$ we use functoriality of Serre
duality in both categories and the definition of $\Phi_\tau$ to reduce it to
the case $k=l=0$. 

Assume $k=l=0$. We first deal with the case where we have
$A_\nu=\mathcal{L}(\varphi_\nu) \otimes F(V_\nu, \exp(N_\nu))$ for all $\nu$.

If $\Lambda_1\ne \Lambda_2\ne \Lambda_3\ne \Lambda_1$, the commutativity of the
diagram was derived in \cite{PZ} from an addition formula for theta functions.

If $\Lambda_1=\Lambda_2\ne \Lambda_3$ we can assume $\varphi_1=\varphi_2$ since
otherwise $\Hom(A_1, A_2)=0$. This implies, using Lemma \ref{flatsect}:
\begin{multline}
\Hom(A_1,A_2) = \Hom(F(V_1,\exp(N_1)), F(V_2,\exp(N_2))) =\notag \\= 
\{f\in\Hom(V_1,V_2)\mid f\circ N_1=N_2\circ f\},  
\end{multline}
which is identified by $\Phi_\tau$ with 
\begin{multline}
\Hom(\Phi(A_1), \Phi(A_2)) = H^0(\Lambda_1, \HOM(M_1, M_2)) =\notag \\ =
\{ f\in \Hom(V_1,V_2)\mid f\circ M_1=M_2\circ f\}.
\end{multline}
On the other hand, if we use the explicit description of the isomorphism
$$\Hom(A_2,A_3) \cong 
H^0(\mathcal{L}(\varphi_2^{-1}\varphi_3^{\phantom{1}}))\otimes\Hom(V_2,V_3)$$
given in \cite[Prop. 2]{PZ}, we see immediately that the composition translates
to the map 
$$
\begin{CD}
  \{f\in\Hom(V_1,V_2)\mid f\circ N_1=N_2\circ f\} \otimes
  H^0(\mathcal{L}(\varphi_2^{-1}\varphi_3^{\phantom{1}})) \otimes \Hom(V_2,V_3)
  \\ 
  @VVV\\
  H^0(\mathcal{AL}(\varphi_1^{-1}\varphi_3^{\phantom{1}})) \otimes
  \Hom(V_1,V_3) 
\end{CD}
$$ sending $f\otimes\theta\otimes g$ to $\theta\otimes(g\circ f)$.
Now, the commutativity
of the above diagram is clear from the definition of $\Phi_\tau$. 
If $\Lambda_1\ne \Lambda_2=\Lambda_3$ the proof is the same, so we skip it. 

If $\Lambda_1=\Lambda_3$ we have seen above that we can have a nonzero
composition only if $\Lambda_1=\Lambda_2=\Lambda_3$ or we can reduce the claim
to the previously discussed situation, using Serre duality. If
$\Lambda_1=\Lambda_2=\Lambda_3$ we can assume $\varphi_1=\varphi_2\varphi_3$
and again from Lemma \ref{flatsect} and the definition of $\Phi_\tau$ the
commutativity is obtained. 

The case where $A_\nu =
\pi_{r_\nu\ast}(\mathcal{L}_{r_\nu\tau}(\varphi_\nu) \otimes
F_{r_\nu\tau}(V_\nu, \exp(N_\nu)))$ are arbitrary indecomposable vector
bundles on $E_\tau$ was studied in \cite{PZ}. Finally, we have to show
commutativity of the diagram above in case at least one of the sheaves $A_\nu$
is an irreducible torsion sheaf. But again by putting together the descriptions
and definitions of the spaces involved the result follows immediately.

We have established that 
$$\Phi_\tau: \mathcal{D}^b(E_\tau) \rightarrow \mathcal{FK}^0(E^\tau)$$ is a
functor. Just by definition this functor is additive, fully faithful and
compatible with shifts. 

To finish the proof of Theorem \ref{equiv} it is
sufficient to show that any indecomposable object in $\mathcal{FK}^0(E^\tau)$
is isomorphic to an object of the form $\Phi_\tau(A)$. Let $(\Lambda,\alpha,M)$
be an indecomposable object in $\mathcal{FK}^0(E^\tau)$. Indecomposability
means that the local system $M$ is indecomposable. Since we assume that the
eigenvalues of the monodromy are of modulus one, we can describe the local
system up to isomorphism by $$M=\exp(-2\pi ib+N)\in\GL(V)$$ with a (cyclic)
nilpotent endomorphism $N$ and a real number $b$. Because $\Phi_\tau$ is
compatible with shifts, we may assume $-\frac{1}{2}<\alpha\le\frac{1}{2}$. If
$\alpha=\frac{1}{2}$ we obtain $(\Lambda,\alpha,M) =
\Phi_\tau(S(-a\tau-b,V,N))$, where $-1< a\le 0$ is the $x$-intercept of a line
in $\mathbb{R}^2$ representing $\Lambda$. If $\alpha<\frac{1}{2}$ we let
$(r,n)$ be the pair of relatively prime nonnegative integers such that $r+in$
is a real multiple of $\exp(i\pi\alpha)$. If we denote by $0\le \frac{ra}{n}<
\frac{1}{n}$ the smallest possible nonnegative $x$-intercept of a line
representing $\Lambda$, we see easily $(\Lambda,\alpha,M) =
\Phi_\tau(\pi_{r\ast}(\mathcal{L}_{r\tau}(\varphi)\otimes
F_{r\tau}(V,\exp(N))))$ 
with $\varphi = t_{ar\tau+b}^\ast \varphi_0^{\phantom{1}}\cdot
\varphi_0^{n-1}$. Observe 
that replacing $a$ by $a+\frac{u}{r}$ yields the same direct image vector
bundle on $E_\tau$, if $0\le u< r$ is an integer. This completes the proof of
Theorem \ref{equiv}.

\bibliographystyle{amsalpha}

\end{document}